%% file: ex_article.tex
\documentclass[siamart,onefignum,onetabnum]{siamart190516}

\input{ex_shared}

\ifpdf
\hypersetup{
  pdftitle={Computational Poisson Geometry I},
  pdfauthor={Evangelista-Alvarado, Su\'arez-Serrato, Ruiz-Pantale\'on}
}
\fi

\begin{document}

\maketitle

\begin{abstract}
 We present a computational toolkit for (local) Poisson-Nijenhuis calculus on manifolds. Our python module \textsf{PoissonGeometry} implements our algorithms, and accompanies this paper. We include two examples of how our methods can be used, one for gauge transformations of Poisson bivectors in dimension 3, and a second one that determines parametric Poisson bivector fields in dimension 4.

\end{abstract}

\begin{keywords}
  Poisson structures, Poisson-Nijenhuis calculus, Symbolic computation, Python.
\end{keywords}

\begin{AMS}
  68W30, 97N80, 53D17
\end{AMS}

    \section{Introduction}
The origin of the concepts in this paper is the analysis of mechanical systems of Sim\'eon Denis Poisson in 1809 \cite{Poisson}. A Poisson manifold is a pair $(M,\Pi)$, with $M$ a smooth manifold and $\Pi$ a contravariant $2$-tensor field (bivector field) on $M$ satisfying the equation
    \begin{equation}\label{EcJacobiPi}
        \cSch{\Pi,\Pi} \,=\, 0,
    \end{equation}
with respect to the Schouten-Nijenhuis bracket $\cSch{\,\, ,\, }$ for multivector fields \cite{Michor-08,Dufour}. Suppose \mbox{$m=\dim{M}$}, fix a local coordinate system \,\mbox{$x=(U;x^{1},\ldots,x^{m})$}\, on $M$. Then $\Pi$ has the following coordinate representation \cite{Lich-77,11}:
    \begin{equation}\label{EcPiCoord}
        \Pi \,=\, \tfrac{1}{2}\,\Pi^{ij}\,\frac{\partial}{\partial{x^{i}}} \wedge \frac{\partial}{\partial{x^{j}}}
            \,=   \sum_{1 \leq i < j \leq m} \Pi^{ij}\,\frac{\partial}{\partial{x^{i}}} \wedge \frac{\partial}{\partial{x^{j}}}
    \end{equation}
Here, the functions \,\mbox{$\Pi^{ij}=\Pi^{ij}(x) \in \Cinf{U}$}\, are called the coefficients of $\Pi$, and \,$\{\partial/\partial{x^{i}}\}$\, is the canonical basis for vector fields on \,\mbox{$U \subseteq M$}.

The Poisson bivector, and its associated bracket, are essential elements in the comprehension of Hamiltonian dynamics \cite{Camille, Dufour}. We recommend interested readers consult the available surveys of this field \cite{WeinsteinS, Kosmann}.

Table \ref{Funs-Algos-Exes} below compiles the functions in our Python module \textsf{PoissonGeometry}\footnote{Our code repository is found at: {\tt https://github.com/appliedgeometry/poissongeometry}.} their corresponding algorithm, and examples where such objects are used in the references. We describe all of our algorithms in section 2. In section 3 we present two applications that illustrate the usefulness of our computational methods. These are, a new result about gauge transformations of Poisson bivector fields in dimension 3 (Proposition \ref{Prop:gauge}), and a description of parametric families of Poisson bivectors in dimension 4 (Lemma \ref{example}).

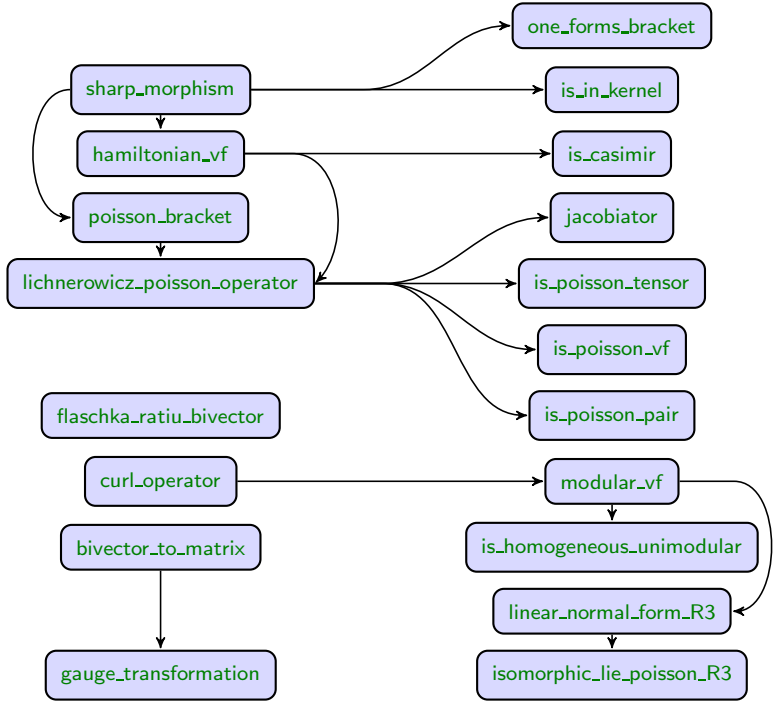
\begin{table}
\centering
{
\renewcommand{\arraystretch}{1.15}
\begin{tabular}{|| l | c | c ||}
     \hline
{\hspace{1.5cm}Function}                   & {Algorithm}                       & {Examples} \\
     \hline
     \hline
{\textsf{sharp\_morphism}}                 & {\ref{AlgSharpMorp}}              & {\cite{Dufour,Camille,Bayro}} \\
     \hline
{\textsf{poisson\_bracket}}                & {\ref{AlgPoissonBracket}}         & {\cite{Camille,Bayro}} \\
     \hline
{\textsf{hamiltonian\_vf}}                 & {\ref{AlgHamVectField}}           & {\cite{Bayro, TV-19}} \\
     \hline
{\textsf{lichnerowicz\_poisson\_operator}} & {\ref{AlgDeltaPi}}                & {\cite{Naka,Poncin}} \\
     \hline
{\textsf{curl\_operator}}                  & {\ref{AlgCurlOperator}}           & {\cite{Damianou,Poncin}} \\
     \hline
{\textsf{bivector\_to\_matrix}}            & {\ref{AlgMatrixPoisson}}          & {\cite{Dufour,Camille,Bayro}} \\
     \hline
{\textsf{jacobiator}}                      & {\ref{AlgJacobiator}}             & {\cite{Dufour,Camille,Bayro}} \\
     \hline
{\textsf{modular\_vf}}                     & {\ref{AlgModularVF}}              & {\cite{Reeb2,4,Poncin}} \\
     \hline
{\textsf{is\_homogeneous\_unimodular}}     & {\ref{AlgHomogUnimod}}            & {\cite{Damianou,Camille,Poncin,Bayro}} \\
     \hline
{\textsf{one\_forms\_bracket}}             & {\ref{AlgOneFormBracket}}         & {\cite{Fer-02, Kosmann}} \\
    \hline
{\textsf{gauge\_transformation}}           & {\ref{AlgGaugeTrans}}             & {\cite{GaugeBursz}} \\
     \hline
{\textsf{linear\_normal\_form\_R3}}        & {\ref{AlgLinNormalFormR3}}        & {\cite{Naka,Bayro}} \\
     \hline
{\textsf{isomorphic\_lie\_poisson\_R3}}    & {\ref{AlgIsomorphicLiePoissonR3}} & {\cite{Naka,Bayro}} \\
     \hline
{\textsf{flaschka\_ratiu\_bivector}}       & {\ref{AlgFlaschka}}               &   {\cite{Damianou, Nar-2015,PabloWrinFib,PabSua-2018}}\\
     \hline
{\textsf{is\_poisson\_tensor}}             & {\ref{AlgIsPoissonTensor}}        & {{\cite{Nar-2015,PabloWrinFib,PabSua-2018}}} \\
     \hline
{\textsf{is\_in\_kernel}}                  & {\ref{AlgIsKernel}}               & {\cite{Dufour,Camille,Poncin,Bayro}} \\
     \hline
{\textsf{is\_casimir}}                     & {\ref{AlgIsCasimir}}              & {{\cite{Damianou,Nar-2015,PabloWrinFib,PabSua-2018}}}\\
     \hline
{\textsf{is\_poisson\_vf}}                 & {\ref{AlgIsPoissonVectorField}}   & {\cite{Naka,MVallYu}} \\
     \hline
{\textsf{is\_poisson\_pair}}               & {\ref{AlgIsPoissonPair}}          & {\cite{MYu,Poncin}} \\
     \hline
\end{tabular}
}
\caption{Functions, corresponding algorithms, and examples where each particular method can be, or has been, applied. The following diagram illustrates functional dependencies in \textsf{PoissonGeometry}.} \label{Funs-Algos-Exes}
\end{table}
\vspace{-1cm}
\begin{center}
\begin{tikzpicture}[
  font=\rmfamily\footnotesize,
  every matrix/.style={ampersand replacement=\&, column sep=2cm, row sep=.20cm},
  source/.style={draw, thick, rounded corners, fill=blue!15, inner sep=.2cm},
  to/.style={->, >=stealth', shorten >=0.5pt, semithick},
  every node/.style={align=center}]

  \matrix{
    {}; \& \node[source] (formsbracket) {\hyperref[AlgOneFormBracket]{\textsf{one\_forms\_bracket}}}; \\
    \node[source] (sharp) {\hyperref[AlgSharpMorp]{\textsf{sharp\_morphism}}};
    \& \node[source] (iskernel) {\hyperref[AlgIsKernel]{\textsf{is\_in\_kernel}}}; \\
    \node[source] (hamiltonian) {\hyperref[AlgHamVectField]{\textsf{hamiltonian\_vf}}};
    \& \node[source] (iscasim) {\hyperref[AlgIsCasimir]{\textsf{is\_casimir}}}; \\
    \node[source] (bracket) {\hyperref[AlgPoissonBracket]{\textsf{poisson\_bracket}}};
    \& \node[source] (jacobiator) {\hyperref[AlgJacobiator]{\textsf{jacobiator}}}; \\
    \node[source] (lichnerowicz) {\hyperref[AlgDeltaPi]{\textsf{lichnerowicz\_poisson\_operator}}};
    \& \node[source] (ispoisson) {\hyperref[AlgIsPoissonTensor]{\textsf{is\_poisson\_tensor}}}; \\
    {}; \& \node[source] (ispoissonvf) {\hyperref[AlgIsPoissonVectorField]{\textsf{is\_poisson\_vf}}}; \\
    \node[source] (flaschka) {\hyperref[AlgFlaschka]{\textsf{flaschka\_ratiu\_bivector}}};
    \& \node[source] (ispoissonpair) {\hyperref[AlgIsPoissonPair]{\textsf{is\_poisson\_pair}}}; \\
    \node[source] (curl) {\hyperref[AlgCurlOperator]{\textsf{curl\_operator}}};
    \& \node[source] (modular) {\hyperref[AlgModularVF]{\textsf{modular\_vf}}}; \\
    \node[source] (matrix) {\hyperref[AlgMatrixPoisson]{\textsf{bivector\_to\_matrix}}};
    \& \node[source] (ishomogunimod) {\hyperref[AlgHomogUnimod]{\textsf{is\_homogeneous\_unimodular}}}; \\
    {}; \& \node[source] (normal) {\hyperref[AlgLinNormalFormR3]{\textsf{linear\_normal\_form\_R3}}}; \\
    \node[source] (gauge) {\hyperref[AlgGaugeTrans]{\textsf{gauge\_transformation}}};
    \& \node[source] (isomorphic) {\hyperref[AlgIsomorphicLiePoissonR3]{\textsf{isomorphic\_lie\_poisson\_R3}}}; \\
        };

  \draw[to] (sharp) -- (hamiltonian);
  \draw[to] (sharp.west) to[out=180, in=180] (bracket.west);
  \draw[to] (sharp.east) --++(0:15mm)to[out=0, in=180] (formsbracket.west);
  \draw[to] (sharp.east) -- (iskernel.west);
  \draw[to] (hamiltonian) --++(0:17.5mm)to[out=0, in=45] (lichnerowicz.east);
  \draw[to] (hamiltonian.east) -- (iscasim.west);
  \draw[to] (bracket) -- (lichnerowicz);
  \draw[to] (lichnerowicz.east) --++(0:9mm)to[out=0, in=180] (jacobiator.west);
  \draw[to] (lichnerowicz.east) -- (ispoisson.west);
  \draw[to] (lichnerowicz.east) --++(0:9mm)to[out=0, in=180] (ispoissonvf.west);
  \draw[to] (lichnerowicz.east) --++(0:9mm)to[out=0, in=180] (ispoissonpair.west);
  \draw[to] (curl.east) -- (modular.west);
  \draw[to] (matrix) -- (gauge);
  \draw[to] (modular) -- (ishomogunimod);
  \draw[to] (modular.east) --++(0:7.23mm)to[out=0, in=0] (normal.east);
  \draw[to] (normal) -- (isomorphic);
\end{tikzpicture}\label{diagram}
\end{center}

    \section{Implementation of Functions in \textsf{PoissonGeometry}} \label{sec:keyf}
In this section we describe the implementation of all functions of the module \textsf{PoissonGeometry}.

    \subsection{Key Functions}
This subsection contains functions that serve as a basis for the implementation of almost all functions of \textsf{PoissonGeometry}.

    \subsubsection{Sharp Morphism} \label{subsec:sharp}
The function \textsf{sharp\_morphism} computes the image of a differential 1-form under the vector bundle morphism \,\mbox{$\Pi^{\natural}: \T^{\ast}M \rightarrow \T{M}$}\, induced by a bivector field $\Pi$ on $M$ and defined by
    \begin{equation}\label{EcPiSharp}
        \big\langle \beta,\Pi^{\natural}(\alpha) \big\rangle \,:=\, \Pi(\alpha,\beta),
    \end{equation}
for any \,\mbox{$\alpha,\beta \in \T^{\ast}M$}\, \cite{Dufour,Camille}. Here, $\langle,\rangle$ is the natural pairing for differential 1-forms and vector fields. Equivalently, \,\mbox{$\Pi^{\natural}(\alpha) = \ii_{\alpha}\Pi$},\, with $\ii_{\bullet}$ the interior product of multivector fields and differential forms defined by the rule \,\mbox{$\ii_{\alpha \wedge \beta} := \ii_{\alpha} \circ \ii_{\beta}$}\, \cite{Kozul}. Analogously for vector fields. In local coordinates, if \,$\alpha = \alpha_{j}\,\dd{x^{j}}$,\, $j=1,\ldots,m$,\, then
    \begin{equation}\label{EcPiSharpCoord}
        \Pi^{\natural}(\alpha) \,= \sum_{1 \leq i<j \leq m} \alpha_{i}\Pi^{ij}\,\frac{\partial}{\partial x^{j}}
                             \,-\, \alpha_{j}\Pi^{ij}\,\frac{\partial}{\partial x^{i}}.
    \end{equation}
\begin{algorithm}[H]
    \captionsetup{justification=centering}
    \caption{\ \textsf{sharp\_morphism}(\emph{bivector, one\_form})} \label{AlgSharpMorp}
        \rule{\textwidth}{0.4pt}
    \Input{a bivector field and a differential 1-form}
    \Output{a vector field which is the image of the differential 1-form under the vector bundle morphism (\ref{EcPiSharp}) induced by the bivector field}
        \rule{\textwidth}{0.4pt}
    \begin{algorithmic}[1] 
        \Procedure{}{}
            \State $m$ $\gets$ dimension of the manifold
                \CommentNew{Given by an instance of \textsf{PoissonGeometry}}
            \State \bluecolor{bivector} $\gets$ a dictionary \{$(1,2)$: $\mathsf{\Pi}^{12}$, ..., $(m-1,m)$: $\mathsf{\Pi}^{m-1\,m}$\} that represents a bivector field according to (\ref{EcMultivectorDic})
            \State \bluecolor{one\_form} $\gets$ a dictionary \{$(1)$: $\mathsf{\salpha}_{1}$, ..., $(m)$: $\mathsf{\salpha}_{m}$\} that represents a differential 1-form according to (\ref{EcMultivectorDic})
            \State \textsc{Convert} each value in \bluecolor{bivector} and in \bluecolor{one\_form} to symbolic expression
            \State \bluecolor{sharp\_dict} $\gets$ the dictionary \{$(1)$: 0,..., $(m)$: 0\}
            \For {each $1 \leq i<j \leq m$}
                    \CommentNew{Compute the sum in (\ref{EcPiSharp})}
                \State \bluecolor{sharp\_dict}[$(i)$] \,$\gets$ \bluecolor{sharp\_dict}[$(i)$]
                    \,$-$\, $\salpha_{j}*\mathsf{\Pi}^{ij}$
                \State \bluecolor{sharp\_dict}[$(j)$] $\gets$ \bluecolor{sharp\_dict}[$(j)$]
                    \,$+$\, $\salpha_{i}*\mathsf{\Pi}^{ij}$
            \EndFor
            \If {all values in \bluecolor{sharp\_dict} are equal to zero}
                \State \textbf{return} \{0: 0\} \CommentNew{A dictionary with zero key and value}
            \Else
                \State \textbf{return} \bluecolor{sharp\_dict}
            \EndIf
        \EndProcedure
    \end{algorithmic}
\end{algorithm}
    \vspace{-0.5cm}
 Observe that the morphism (\ref{EcPiSharp}) is defined, in particular, for Poisson bivector fields. So the function \textsf{sharp\_morphism} can be applied on this class of bivector fields.

    \subsubsection{Poisson Brackets} \label{subsec:bracket}
A Poisson bracket on $M$ is a Lie bracket structure $\{,\}$ on the space of smooth functions $\Cinf{M}$ which is compatible with the pointwise product by the Leibniz rule \cite{Dufour,Camille}. Explicitly, the Poisson bracket induced by a Poisson bivector field $\Pi$ on $M$ is given by the formula
    \begin{equation}\label{EcPiBracket}
        \{f,g\}_{\Pi} \,=\, \big\langle \dd{g},\Pi^{\natural}(\dd{f}) \big\rangle
                      \,=\, \left(\Pi^{\natural}\dd{f}\right)^{i}\frac{\partial{g}}{\partial{x^{i}}}, \qquad \forall\, f,g \in \Cinf{M};
    \end{equation}
for \,\mbox{$i=1,\ldots,m$}. The function \textsf{poisson\_bracket} computes the poisson bracket, induced by a Poisson bivector field, of two scalar functions.
\begin{algorithm}[H]
    \captionsetup{justification=centering}
    \caption{\ \textsf{poisson\_bracket}(\emph{bivector, function\_1, function\_2})} \label{AlgPoissonBracket}
        \rule{\textwidth}{0.4pt}
    \Input{a Poisson bivector field and two scalar functions}
    \Output{the Poisson bracket of the two scalar functions induced by the Poisson bivector field}
        \rule{\textwidth}{0.4pt}
    \begin{algorithmic}[1] 
        \Procedure{}{}
            \State $m$ $\gets$ dimension of the manifold
                \CommentNew{Given by an instance of \textsf{PoissonGeometry}}
            \State \bluecolor{bivector} $\gets$ a dictionary that represents a bivector field according to (\ref{EcMultivectorDic})
            \State \bluecolor{function\_1}, \bluecolor{function\_2} $\gets$ string expressions
                \CommentNew{Represents scalar functions}
            \If {\bluecolor{function\_1} $==$ \bluecolor{function\_2}}
                \State \textbf{return} 0
                    \CommentNew{If $f=g$ in (\ref{EcPiBracket}), then $\{f,g\}=0$}
            \Else
                \State \textsc{Convert} \bluecolor{function\_1} and \bluecolor{function\_2} to symbolic expressions
                \State \bluecolor{gradient\_function\_1} $\gets$ a dictionary that represents the gradient vector field of \bluecolor{ function\_1} according to (\ref{EcMultivectorDic})
                \State \bluecolor{sharp\_function\_1} $\gets$ \textsf{sharp\_morphism}(\bluecolor{bivector}, \bluecolor{ gradient\_function\_1})
                \Statex \CommentNew{See Algorithm \ref{AlgSharpMorp}}
                \State \bluecolor{bracket} $\gets$ 0
                \For {$i = 1 \ \text{to} \ m$}
                    \State \bluecolor{bracket} $\gets$ \bluecolor{bracket}
                        $+$ \bluecolor{sharp\_function\_1}[$(i)$] * {$\partial$(\bluecolor{function\_2})}/{$\partial$\textsf{x}$i$}
                    \Statex \CommentNew{See (\ref{EcPiBracket})}
                \EndFor
                \State \textbf{return} \bluecolor{bracket}
            \EndIf
        \EndProcedure
    \end{algorithmic}
\end{algorithm}
\vspace{-.5cm}

    \subsubsection{Hamiltonian Vector Fields} \label{subsec:hamiltonian}
The function \textsf{hamiltonian\_vf} computes the Hamiltonian vector field
\begin{equation}\label{EcPiHamVectField}
        X_{h} \,:=\, \Pi^{\natural}(\dd{h}).
    \end{equation}
of a function \,\mbox{$h \in \Cinf{M}$}\, respect to a Poisson bivector field $\Pi$ on $M$ \cite{Dufour,Camille}.
\begin{algorithm}[H]
    \captionsetup{justification=centering}
    \caption{\ \textsf{hamiltonian\_vf}(\emph{bivector, hamiltonian\_function})} \label{AlgHamVectField}
        \rule{\textwidth}{0.4pt}
    \Input{a Poisson bivector field and a scalar function}
    \Output{the Hamiltonian vector field of the scalar function relative to the Poisson bivector field}
        \rule{\textwidth}{0.4pt}
    \begin{algorithmic}[1] 
        \Procedure{}{}
            \State $m$ $\gets$ dimension of the manifold
                \CommentNew{Given by an instance of \textsf{PoissonGeometry}}
            \State \bluecolor{bivector} $\gets$ a dictionary that represents a bivector field according to (\ref{EcMultivectorDic})
            \State \bluecolor{hamiltonian\_function} $\gets$ a string expression
                \CommentNew{Represents a scalar function}
                    \breakalg
            \State \textsc{Convert} \bluecolor{hamiltonian\_function} to symbolic expression
                    \State \bluecolor{gradient\_hamiltonian} $\gets$ a dictionary that represents the gradient vector field of \bluecolor{ hamiltonian\_function} according to (\ref{EcMultivectorDic})
                    \State \textbf{return} \textsf{sharp\_morphism}(\bluecolor{bivector}, \bluecolor{gradient\_hamiltonian})
                    \Statex \CommentNew{See Algorithm \ref{AlgSharpMorp} and formula (\ref{EcPiHamVectField})}
        \EndProcedure
    \end{algorithmic}
\end{algorithm}
\vspace{-.5cm}

    \subsubsection{Coboundary Operator} \label{subsec:coboundary}
The adjoint operator of a Poisson bivector field $\Pi$ on $M$ with respect to the Schouten-Nijenhuis bracket gives rise to a cochain complex \mbox{$(\Gamma\wedge\T{M},\delta_{\pi})$}, called the \emph{Lichnerowicz-Poisson} complex of $(M,\Pi)$ \cite{Lich-77,Dufour,Camille}.  Here \,$\delta_{\Pi}:\Gamma\wedge^{\bullet}\T{M} \rightarrow \Gamma\wedge^{\bullet+1}\T{M}$\, is the coboundary operator ($\delta_{\Pi}^{2}=0$) defined by
     \begin{equation}\label{EcDeltaPi}
          \delta_{\Pi}(A) \,:=\, \cSch{\Pi,A}, \qquad \forall\, A \in \Gamma\wedge\T{M}.
     \end{equation}
Here, \mbox{$\Gamma\wedge\T{M}$} denotes the $\Cinf{M}$--module of multivector fields on $M$. Explicitly, if \,$a=\deg{A}$,\, then for any \,\mbox{$f_{1},\ldots,f_{a+1} \in \Cinf{M}$}\,:
    \begin{multline*}
        \cSch{\Pi,A}(\dd{f_{1}},\ldots,\dd{f_{a+1}}) \,=\, \sum_{k=1}^{a+1}(-1)^{k+1}\, \big\{f_{k},A(\dd{f_{1}},\ldots,\widehat{\dd{f_{k}}},\ldots,\dd{f_{a+1}})\big\}_{\Pi} \\
            + \sum_{1 \leq k < l \leq a+1}(-1)^{k+l}\, A\big(\dd{\{f_{k},f_{l}\}}_{\Pi},\dd{f_{1}},\ldots,\widehat{\dd{f_{k}}},\ldots,\widehat{\dd{f_{l}}},\ldots,\dd{f_{a+1}}\big)
    \end{multline*}
Throughout this paper the symbol \,$\widehat{\ }$\, will denote the absence of the corresponding factor.  In particular, if \,$f_{1}=x^{1},\ldots,f_{a+1}=x^{a+1}$\, are local coordinates on $M$, we have that for \,$1 \leq i_{1} < \cdots< i_{a+1} \leq m$\,:
    \begin{align}\label{EcDeltaPiCoord}
        \cSch{\Pi,A}^{i_{1} \cdots i_{a+1}} \,&= \sum_{k=1}^{a+1}(-1)^{k+1}\,\big\{x^{i_{k}},A^{i_{1} \cdots \widehat{i_{k}} \cdots i_{a+1}}\big\}_{\Pi} \\ 
            &+ \hspace{-0.11cm} \sum_{1 \leq k < l \leq a+1}(-1)^{k+l}\, \frac{\partial \Pi^{i_{k}i_{l}}}{\partial x^{s}}\,A^{s\,i_{1} \cdots \widehat{i_{k}} \cdots \widehat{i_{l}} \cdots i_{a+1}}
    \end{align}
Here $\cSch{\Pi,A}^{i_{1} \cdots i_{a+1}}$, $A^{i_{1} \cdots i_{a}}$ and $\Pi^{i_{k}i_{l}}$ are the coefficients of the coordinate expressions of $\cSch{\Pi,A}$, $A$ and $\Pi$, in that order. The function \textsf{lichnerowicz\_poisson\_operator} computes the image of a multivector field under the coboundary operator induced by a Poisson bivector field.
\begin{algorithm}[H]
    \captionsetup{justification=centering}
    \caption{\ \textsf{lichnerowicz\_poisson\_operator}(\emph{bivector, multivector})} \label{AlgDeltaPi}
        \rule{\textwidth}{0.4pt}
    \Input{a Poisson bivector field and a multivector field}
    \Output{the image of the multivector field under the coboundary operator (\ref{EcDeltaPi}) induced by the Poisson bivector field}
        \rule{\textwidth}{0.4pt}
    \begin{algorithmic}[1] 
        \Procedure{}{}
            \State $m$ $\gets$ dimension of the manifold
                \CommentNew{Given by an instance of \textsf{PoissonGeometry}}
            \State \bluecolor{a} $\gets$ degree of the multivector field
            \State \bluecolor{bivector} $\gets$ a dictionary \{$(1,2)$: $\mathsf{\Pi}^{12}$, ..., $(m-1,m)$: $\mathsf{\Pi}^{m-1\,m}$\} that represents a bivector field according to (\ref{EcMultivectorDic})
            \State \bluecolor{multivector} $\gets$ a dictionary \mbox{\{($1,...,a$): $\mathsf{A}^{1 \cdots a}$, ..., ($m$-$a+1,...,m$): $\mathsf{A}^{m-a+1 \cdots m}$\}} that represents a bivector field according to (\ref{EcMultivectorDic}) or a string expression
                \breakalg
            \If {\bluecolor{a} $+\, 1 > m$}
                    \CommentNew{$\deg\,\cSch{\Pi,A} = \deg(A) + 1$ in (\ref{EcDeltaPi})}
                \State \textbf{return} \{0: 0\}
                    \CommentNew{A dictionary with zero key and value}
            \ElsIf {\bluecolor{multivector} is a string expression}
                \State \textbf{return} \textsf{hamiltonian\_vf}(\bluecolor{bivector}, \textsf{str}($-1$*) \textsf{+} \bluecolor{multivector})
                \Statex \CommentNew{See Algorithm \ref{AlgHamVectField}}
            \Else
                \State \textsc{Convert} each value in \bluecolor{bivector} and in \bluecolor{multivector} to symbolic expression
                \State \bluecolor{image\_multivector} $\gets$ the dictionary \mbox{\{($1,...,a+1$): 0, ..., ($m-a,...,m$): 0\}}
                \For {each $1 \leq i_{1} < \cdots < i_{\text{\bluecolor{a}} + 1} \leq m$}
                    \For {$k = 1$ to \bluecolor{a}\,$+1$} \CommentNew{Compute first summation in (\ref{EcDeltaPiCoord})}
                        \State \bluecolor{image\_multivector}[$(i_{1},..., i_{\text{\bluecolor{a}} + 1})$] $\gets$ \bluecolor{image\_multivector}[$(i_{1},..., i_{\text{\bluecolor{a}} + 1})$]
                            $+$ $(-1)^{k+1}\,\big\{\textsf{x}{i_{k}},\mathsf{A}^{i_{1} \cdots \widehat{i_{k}} \cdots i_{a+1}}\big\}_{\Pi}$
                                \CommentNew{See Algorithm \ref{AlgPoissonBracket} to compute the Poisson bracket $\{,\}_{\Pi}$}
                    \EndFor
                    \For {each $1 \leq k < l \leq $ \bluecolor{a}\,$+1$}
                            \CommentNew{Compute second summation in (\ref{EcDeltaPiCoord})}
                        \State \bluecolor{image\_multivector}[$(i_{1},..., i_{\text{\bluecolor{a}} + 1})$] $\gets$ \bluecolor{ image\_multivector}[$(i_{1},..., i_{\text{\bluecolor{a}} + 1})$]
                        \,$+$\, $(-1)^{k+l}\,\frac{\partial \mathsf{\Pi}^{i_{k}i_{l}}}{\partial \textsf{x}{s}}\,\mathsf{A}^{s\,i_{1} \cdots \widehat{i_{k}} \cdots \widehat{i_{l}} \cdots i_{a+1}}$
                    \EndFor
                \EndFor
                \If {all values in \bluecolor{image\_multivector} are equal to zero}
                    \State \textbf{return} \{0: 0\} \CommentNew{A dictionary with zero key and value}
                \Else
                    \State \textbf{return} \bluecolor{image\_multivector}
                \EndIf
            \EndIf
        \EndProcedure
    \end{algorithmic}
\end{algorithm}
\vspace{-.5cm}

    \subsubsection{Curl (Divergence) Operator} \label{subsec:curl}
Fix a volume form $\Omega_{0}$ on an oriented Poisson manifold $(M, \Pi,$ $\Omega_{0})$. The \emph{divergence} (relative to $\Omega_{0}$) of an $a$--multivector field $A$ on $M$ is the unique $(a-1)$--multivector field $\mathscr{D}_{\Omega_{0}}(A)$ on $M$ such that
    \begin{equation}\label{EcTrazaDef}
        \ii_{\mathscr{D}_{\Omega_{0}}(A)}\Omega_{0} \,=\, \dd{\ii_{A}}\Omega_{0}.
    \end{equation}
This induces a (well defined, $\Omega_{0}$--dependent) coboundary operator \,\mbox{$\mathscr{D}_{\Omega_{0}}:A \mapsto \mathscr{D}_{\Omega_{0}}(A)$}\, on the module of multivector fields on $M$, called the \emph{curl operator} \cite{Kozul,Camille}. As any other volume form on $M$ is a multiple $f\Omega_{0}$ of $\Omega_{0}$ by a nowhere vanishing function \mbox{$f \in \Cinf{M}$}, we have \,\mbox{$\mathscr{D}_{f\Omega_{0}}=\mathscr{D}_{\Omega_{0}}+\tfrac{1}{f}\ii_{\dd{f}}$}.\,  In local coordinates, expressing $\Omega_{0}$ as
    \begin{equation}\label{EcOmegaCero}
        \Omega_{0} \,=\, \dd{x^{1}} \wedge \cdots \wedge \dd{x^{m}},
    \end{equation}
 then for any $a$--multivector field \,\mbox{$A=A^{i_{1} \cdots i_{a}}\partial/\partial{x^{i_{1}}} \wedge \cdots \wedge \partial/\partial{x^{i_{a}}}$}\, on $M$, with \,{$1 \leq i_{1} < \cdots< i_{a} \leq m$},\, the divergence of $A$ with respect to the volume form $f\Omega_{0}$ is given by:
    \begin{equation}\label{EcTraza}
        \mathscr{D}_{f\Omega_{0}}(A) = \sum_{k=1}^{m} (-1)^{k + 1}\left(\frac{\partial{A^{i_{1} \cdots i_{a}}}}{\partial{x^{i_{k}}}}
            + \frac{1}{f}\frac{\partial{f}}{\partial{x^{i_{k}}}}\,A^{i_{1} \cdots i_{a}}\right) \frac{\partial}{\partial{x^{i_{1}}}} \wedge \cdots \wedge \widehat{\frac{\partial}{\partial{x^{i_{k}}}}} \wedge \cdots \wedge \frac{\partial}{\partial{x^{i_{a}}}}
    \end{equation}

Let  $f_{0}$ be a nonzero scalar function. The function \textsf{curl\_operator} computes the divergence of a multivector field respect to the volume form $f_{0}\Omega_{0}$, for $\Omega_{0}$ in (\ref{EcOmegaCero}).
\begin{algorithm}[H]
    \captionsetup{justification=centering}
    \caption{\ \textsf{curl\_operator}(\emph{multivector, function})} \label{AlgCurlOperator}
        \rule{\textwidth}{0.4pt}
    \Input{a multivector field and a nonzero scalar function $f_{0}$}
    \Output{the divergence of the multivector field with respect to the volume form $f_{0}\Omega_{0}$}
        \rule{\textwidth}{0.4pt}
    \begin{algorithmic}[1] 
        \Procedure{}{}
            \State $m$ $\gets$ dimension of the manifold
                \CommentNew{Given by an instance of \textsf{PoissonGeometry}}
            \State \bluecolor{a} $\gets$ degree of the multivector field
            \State \bluecolor{multivector} $\gets$ a dictionary \mbox{\{($1$,...,$a$): $\mathsf{A}^{1 \cdots a}$, ..., ($m$-$a+1$,...,$m$): $\mathsf{A}^{m-a+1 \cdots m}$\}} that represents a multivector field according to (\ref{EcMultivectorDic}) or a string expression
            \State \bluecolor{function} $\gets$ a string expression
                \CommentNew{Represents a nonzero function}
            \If {\bluecolor{multivector} is a string expression}
                \State \textbf{return} \{0: 0\} \CommentNew{A dictionary with zero key and value}
            \Else
                \State \textsc{Convert} each value in \bluecolor{multivector} and \bluecolor{function} to symbolic expression
                \State \bluecolor{curl\_multivec} $\gets$ the dictionary \mbox{\{($1$,...,$a-1$): 0, ..., ($m$-$a+2$,...,$m$): 0\}}
                \For {each $1 \leq i_{1} < \cdots < i_{\text{\bluecolor{a}}} \leq m$}
                        \CommentNew{Compute the summation in (\ref{EcTraza})}
                    \For {$k=1$ to $m$}
                        \State \bluecolor{curl\_multivec}[$(i_{1},...,\widehat{\imath_{k}},...,i_{\text{\bluecolor{a}}})$] $\gets$ \bluecolor{curl\_multivec}[$(i_{1},...,\widehat{\imath_{k}},...,i_{\text{\bluecolor{a}}})$]
                        + $(-1)^{k+1}$ * \Big($\frac{\partial{\text{\textsf{A}}^{i_{1} \cdots i_{a}}}}{\partial{\text{\textsf{x}}{i_{k}}}}$ \,+\, $\frac{1}{\text{\bluecolor{function}}}$ * $\frac{\partial{(\text{\bluecolor{ function}})}}{\partial{\text{\textsf{x}}{i_{k}}}}$ * \text{\textsf{A}}$^{i_{1} \cdots i_{a}}$\Big)
                     \EndFor
                \EndFor
                \If {all values in \bluecolor{curl\_multivec} are equal to zero}
                    \State \textbf{return} \{0: 0\} \CommentNew{A dictionary with zero key and value}
                \Else
                    \State \textbf{return} \bluecolor{curl\_multivec}
                \EndIf
            \EndIf
        \EndProcedure
    \end{algorithmic}
\end{algorithm}
\vspace{-.5cm}

  \subsection{Matrix of a bivector field} \label{subsec:matrix}
The function \textsf{bivector\_to\_matrix} computes the (local) matrix $\left[\Pi^{ij}\right]$ of a bivector field $\Pi$ on $M$ \cite{Dufour,Camille}, the coefficients of $\Pi$ in (\ref{EcPiCoord}). In particular, it computes the matrix of a Poisson bivector field.
\begin{algorithm}[H]
    \captionsetup{justification=centering}
    \caption{\ \textsf{bivector\_to\_matrix}(\emph{bivector})} \label{AlgMatrixPoisson}
        \rule{\textwidth}{0.4pt}
    \Input{a bivector field}
    \Output{the (local) matrix of the bivector field} \\[-0.25cm]
        \rule{\textwidth}{0.4pt}
\begin{algorithmic}[1] 
    \Procedure{}{}
            \State $m$ $\gets$ dimension of the manifold
                \CommentNew{Given by an instance of \textsf{PoissonGeometry}}
            \State \bluecolor{bivector} $\gets$ a dictionary \{$(1,2)$: $\mathsf{\Pi}^{12}$, ..., $(m-1,m)$: $\mathsf{\Pi}^{m-1\,m}$\} that represents a bivector field according to (\ref{EcMultivectorDic})
            \State \bluecolor{matrix} $\gets$ a symbolic $m \times m$-matrix
        \For {each $1 \leq i<j \leq m$}
            \State \textsc{Convert} $\mathsf{\Pi}^{ij}$ to symbolic expression
            \State \bluecolor{matrix}[$i-1,j-1$] $\gets$ $\mathsf{\Pi}^{ij}$
            \State \bluecolor{matrix}[$j-1,i-1$] $\gets$ $(-1)$ * \bluecolor{matrix}[$i-1,j-1$]
        \EndFor
            \breakalg
        \State \textbf{return} \bluecolor{matrix}
    \EndProcedure
    \end{algorithmic}
\end{algorithm}
\vspace{-.5cm}

    \subsection{Jacobiator} \label{subsec:jacobiator}
The Schouten-Nijenhuis bracket of a bivector field $\Pi$ with itself, $\cSch{\Pi,\Pi}$, is computed with the \textsf{jacobiator} function. This 3-multivector field is called the Jacobiator of $\Pi$. The Jacobi identity (\ref{EcJacobiPi}) for $\Pi$ follows from the vanishing of its Jacobiator \cite{Dufour,Camille}.
\begin{algorithm}[H]
    \captionsetup{justification=centering}
    \caption{\ \textsf{jacobiator}(\emph{bivector})} \label{AlgJacobiator}
        \rule{\textwidth}{0.4pt}
    \Input{a bivector field}
    \textbf{Output:} {the Schouten-Nijenhuis bracket of the bivector field with itself} \\[-0.25cm]
        \rule{\textwidth}{0.4pt}
    \begin{algorithmic}[1] 
    \Procedure{}{}
            \State \bluecolor{bivector} $\gets$ a dictionary that represents a bivector field according to (\ref{EcMultivectorDic})
            \State \textbf{return} \textsf{lichnerowicz\_poisson\_operator}(\bluecolor{bivector}, \bluecolor{bivector}) \CommentNew{See Algorithm \ref{AlgDeltaPi}}
\EndProcedure
    \end{algorithmic}
\end{algorithm}
\vspace{-.5cm}

    \subsection{Modular Vector Field} \label{subsec:modular}
For $(M,\Pi,\Omega)$ an orientable Poisson manifold, and a fixed volume form $\Omega$ on $M$, the map
     \begin{equation}\label{EcModularVF}
          Z: h \,\longmapsto\, \mathscr{D}_{\Omega}(X_{h})
     \end{equation}
is a derivation of $\Cinf{M}$. Therefore it defines a vector field on $M$, called the \emph{modular vector field of $\Pi$} relative to $\Omega$ \cite{WeModular,Reeb2,Dufour,Camille}. Here, $\mathscr{D}_{\Omega}$ is the curl operator relative to $\Omega$ (\ref{EcTrazaDef}). Then, $Z$ is a Poisson vector field of $\Pi$ which is independent of the choice of a volume form, {\em modulo} Hamiltonian vector fields: \,\mbox{$Z_{f\Omega}=Z - \tfrac{1}{f}X_{f}$}.\, Here $Z_{f\Omega}$ is the modular vector field of $\Pi$ relative to the volume form $f\Omega$ and \,\mbox{$f \in \Cinf{M}$}\, a nowhere vanishing function. In this context, the Poisson bivector field $\Pi$ is said to be \emph{unimodular} if $Z$ is a Hamiltonian vector field (\ref{EcPiHamVectField}). Equivalently, if $Z$ is zero for some volume form on $M$.
We can compute the modular vector field $Z$ of $\Pi$ (\ref{EcModularVF}) relative to a volume form $f\Omega$ as the (minus) divergence of $\Pi$ (\ref{EcTraza}):
    \begin{equation}\label{EcModuVFTraza}
        Z_{f\Omega} \,=\, -\mathscr{D}_{f\Omega}(\Pi) \,=\, \mathscr{D}_{f\Omega}(-\Pi).
    \end{equation}
Let $f_{0}$ be a nonzero scalar function. The function \textsf{modular\_vectorfield}  computes the modular vector field of a Poisson bivector field with respect to the volume form $f_{0}\Omega_{0}$.
\begin{algorithm}[H]
    \captionsetup{justification=centering}
    \caption{\ \textsf{modular\_vf}(\emph{bivector, function})} \label{AlgModularVF}
        \rule{\textwidth}{0.4pt}
    \Input{a Poisson bivector field and a nonzero scalar function $f_{0}$}
    \Output{the modular vector field of the Poisson bivector field (\ref{EcModuVFTraza}) relative to the volume form $f_{0}\Omega_{0}$, }
        \rule{\textwidth}{0.4pt}
    \begin{algorithmic}[1] 
        \Procedure{}{}
            \State \bluecolor{bivector} $\gets$ a dictionary \{$(1,2)$: $\mathsf{\Pi}^{12}$, ..., $(m-1,m)$: $\mathsf{\Pi}^{m-1\,m}$\} that represents a bivector field according to (\ref{EcMultivectorDic})
            \State \bluecolor{function} $\gets$ a string expression
                \CommentNew{Represents a scalar function}
            \For{each $1 \leq i < j \leq m$}
                    \CommentNew{A dictionary for $-\Pi$ in (\ref{EcModuVFTraza})}
                \State \bluecolor{bivector}[($i,j$)] $\gets$ $-\,\mathsf{\Pi}^{ij}$
            \EndFor
                \breakalg
            \State \textbf{return} \textsf{curl\_operator}(\bluecolor{bivector}, \bluecolor{function})
                \CommentNew{See Algorithm \ref{AlgCurlOperator} and formula (\ref{EcModuVFTraza})}
        \EndProcedure
    \end{algorithmic}
\end{algorithm}
\vspace{-.5cm}

    \subsection{Unimodularity of Homogeneous Poisson bivector fields} \label{subsec:unimodhomog}
We can verify whether an homogeneous Poisson bivector field is unimodular or not with the \textsf{is\_homogeneous\_unimodular} function.  A Poisson bivector field $\Pi$ on $\R{m}_{x}$,
    \begin{equation}\label{EcPiHomog}
          \Pi \,=\, \tfrac{1}{2}\Pi^{ij}\,\frac{\partial}{\partial{x_{i}}} \wedge \frac{\partial}{\partial{x_{j}}}, \qquad i,j=1\ldots,m;
     \end{equation}
is said to be \emph{homogeneous} if each coefficient $\Pi^{ij}$ is an homogeneous polynomial \cite{Camille}. To implement this function we use the following fact: an homogeneous Poisson bivector field on $\R{m}_{x}$ is unimodular on (the whole of) $\R{m}_{x}$ if and only if its modular vector field (\ref{EcModularVF}) relative to the Euclidean volume form is zero \cite{Koz-95}.
\begin{algorithm}[H]
    \captionsetup{justification=centering}
    \caption{\ \textsf{is\_homogeneous\_unimodular}(\emph{bivector})} \label{IsHomogUnimod} \label{AlgHomogUnimod}
        \rule{\textwidth}{0.4pt}
    \textbf{Input:} a homogeneous Poisson bivector field on $\R{m}$ \\
    \Output{verify if the modular vector field respect to the Euclidean volume form on $\R{m}$ of the Poisson bivector field is zero or not}
        \rule{\textwidth}{0.4pt}
\begin{algorithmic}[1] 
    \Procedure{}{}
            \State \bluecolor{bivector} $\,\gets\,$ a dictionary that represents a bivector field according to (\ref{EcMultivectorDic})
            \If {\textsf{modular\_vf}(\bluecolor{bivector}, 1) == \textsf{\{0: 0\}}}
                    \CommentNew{See Algorithm \ref{AlgModularVF}}
                \State \textbf{return} True
            \Else
                \State \textbf{return} False
            \EndIf
\EndProcedure
    \end{algorithmic}
\end{algorithm}
\vspace{-.5cm}

    \subsection{Bracket on Differential 1-Forms} \label{subsec:oneforbracket}
The function \textsf{one\_forms\_bracket} computes the Lie bracket of two differential $1$-forms \,$\alpha, \beta \in \Gamma\,\T^{\ast}M$\, induced by a Poisson bivector field $\Pi$ on $M$ \cite{Dufour,Camille} and defined by
    \begin{equation*}
        \{\alpha,\beta\}_{\Pi} \,:=\, \ii_{\Pi^{\natural}(\alpha)}\dd{\beta} - \ii_{\Pi^{\natural}(\beta)}\dd{\alpha} + \dd\big\langle \beta,\Pi^{\natural}(\alpha) \big\rangle.
    \end{equation*}
Here, $\dd$ is the exterior derivative for differential forms and \,\mbox{$\{\dd{f},\dd{g}\}_{\Pi} = \dd\{f,g\}_{\Pi}$},\, by definition for all \,\mbox{$f,g \in \Cinf{M}$}.\, The bracket on the right-hand side of this equality is the Poisson bracket for smooth functions on $M$ induced by $\Pi$ (\ref{EcPiBracket}). In coordinates, if \,$\alpha = \alpha_{k}\,\dd{x^{k}}$\, and \,$\beta = \beta_{l}\,\dd{x^{l}}$,\, for \,$k,l=1,\ldots,m$\,:
    \begin{multline}\label{EcOneFormBracket}
        \{\alpha,\beta\}_{\Pi} \,= \sum_{1 \leq i < j \leq m}\left[ \big( \Pi^{\natural}\alpha \big)^{j}\left( \frac{\partial{\beta_{i}}}{\partial{x^j}}
            - \frac{\partial{\beta_{j}}}{\partial{x^i}} \right)
            - \big( \Pi^{\natural}\beta \big)^{j}\left( \frac{\partial{\alpha_{i}}}{\partial{x^j}}
            - \frac{\partial{\alpha_{j}}}{\partial{x^i}} \right) \right] \dd{x^{i}} \\
            +\sum_{1 \leq i < j \leq m}\left[ \big( \Pi^{\natural}\alpha \big)^{i}\left( \frac{\partial{\beta_{j}}}{\partial{x^i}}
            - \frac{\partial{\beta_{i}}}{\partial{x^j}} \right)
            - \big( \Pi^{\natural}\beta \big)^{i}\left( \frac{\partial{\alpha_{j}}}{\partial{x^i}}
            - \frac{\partial{\alpha_{i}}}{\partial{x^j}} \right) \right]\dd{x^{j}}
            + \frac{\partial\left[(\Pi^{\natural}\alpha)^{l}\beta_{l}\right]}{\partial{x^{k}}}\dd{x^{k}}
    \end{multline}
\begin{algorithm}[H]
    \captionsetup{justification=centering}
    \caption{\ \textsf{one\_forms\_bracket}(\emph{bivector, one\_form\_1, one\_form\_2})} \label{AlgOneFormBracket}
        \rule{\textwidth}{0.4pt}
    \Input{a Poisson bivector field and two differential 1-forms}
    \Output{a differential 1-form which is the Lie bracket induced by the Poisson bivector field of the two differential 1-forms}
        \rule{\textwidth}{0.4pt}
    \begin{algorithmic}[1] 
        \Procedure{}{}
            \State $m$ $\gets$ dimension of the manifold
                \CommentNew{Given by an instance of \textsf{PoissonGeometry}}
            \State \bluecolor{bivector} $\gets$ a dictionary that represents a bivector field according to (\ref{EcMultivectorDic})
            \State \bluecolor{one\_form\_1} $\gets$ a dictionary \{$(1)$: $\mathsf{\salpha}_{1}$, ..., $(m)$: $\mathsf{\salpha}_{m}$\} that represents a differential 1-form according to (\ref{EcMultivectorDic})
            \State \bluecolor{one\_form\_2} $\gets$ a dictionary \{$(1)$: $\mathsf{\sbeta}_{1}$, ..., $(m)$: $\mathsf{\beta}_{m}$\} that represents a differential 1-form according to (\ref{EcMultivectorDic})
            \State \bluecolor{sharp\_1} $\gets$ \textsf{sharp\_morphism}(\bluecolor{bivector}, \bluecolor{one\_form\_1})
                \CommentNew{See Algorithm \ref{AlgSharpMorp}}
            \State \bluecolor{sharp\_2} $\gets$ \textsf{sharp\_morphism}(\bluecolor{bivector}, \bluecolor{one\_form\_2})
            \State \textsc{Convert} each of \bluecolor{one\_form\_1} and  \bluecolor{one\_form\_2} to a symbolic expression
            \State \bluecolor{forms\_bracket} $\gets$ the dictionary \{$(1)$: 0,..., $(m)$: 0\}
            \For{$1 \leq i<j \leq m$}
                    \CommentNew{Compute the first two summations in (\ref{EcOneFormBracket})}
                \State \bluecolor{forms\_bracket}[$(i)$] $\gets$ \bluecolor{forms\_bracket}[$(i)$]
                    + \bluecolor{ sharp\_1}[$(j)$] * $\big({\partial\mathsf{\sbeta}_{i}}/{\partial\text{\textsf{x}}j}
                    - {\partial\mathsf{\sbeta}_{j}}/{\partial\text{\textsf{x}}i}\big)$ $-$ \bluecolor{sharp\_2}[$(j)$]
                    * $\big({\partial\mathsf{\salpha}_{i}}/{\partial\text{\textsf{x}}j}
                    - {\partial\mathsf{\salpha}_{j}}/{\partial\text{\textsf{x}}i}\big)$
                \State \bluecolor{forms\_bracket}[$(j)$] $\gets$ \bluecolor{forms\_bracket}[$(j)$]
                    + \bluecolor{ sharp\_1}[$(i)$] * $\big({\partial\mathsf{\sbeta}_{j}}/{\partial\text{\textsf{x}}i}
                    - {\partial\mathsf{\sbeta}_{i}}/{\partial\text{\textsf{x}}j}\big)$ $-$ \bluecolor{sharp\_2}[$(i)$]
                    * $\big({\partial\mathsf{\salpha}_{j}}/{\partial\text{\textsf{x}}i}
                    - {\partial\mathsf{\salpha}_{i}}/{\partial\text{\textsf{x}}j}\big)$
            \EndFor
            \For{$k,l=1$ to $m$}
                    \CommentNew{Compute the last sum in (\ref{EcOneFormBracket})}
                \State \bluecolor{forms\_bracket}[$(k)$] $\gets$ \bluecolor{forms\_bracket}[$(k)$]
                    + {$\partial$(\bluecolor{ sharp\_1}[$(l)$] * $\mathsf{\sbeta}_{l}$)}/{$\partial$\textsf{x}$k$}
            \EndFor
            \If {all values in \bluecolor{forms\_bracket} are equal to zero}
                \State \textbf{return} \{0: 0\} \CommentNew{A dictionary with zero key and value}
            \Else
                \State \textbf{return} \bluecolor{forms\_bracket}
            \EndIf
        \EndProcedure
    \end{algorithmic}
\end{algorithm}
\vspace{-.5cm}

    \subsection{Gauge Transformations} \label{subsec:gauge}
Let $\Pi$ be a bivector field on $M$. Suppose we are given a differential 2-form $\lambda$ on $M$ such that the vector bundle morphism
    \begin{equation}\label{EcGaugeInvertible}
          \big(\mathrm{id}_{\T^{\ast}{M}} - \lambda^{\flat} \circ \Pi^{\natural}\big):\T^{\ast}M \rightarrow \T{M} \qquad \text{\rm{is \ invertible}}.
     \end{equation}
Then, there exists a bivector field $\overline{\Pi}$ on $M$ (well) defined by the skew-symmetric morphism
    \begin{equation}\label{EcGaugeEquiv}
          \overline{\Pi}^{\natural} \,=\, \Pi^{\natural} \circ \big(\, \mathrm{id}_{\T^{\ast}{M}}-\lambda^{\flat} \circ \Pi^{\natural} \,\big)^{-1}.
     \end{equation}
Here, \,\mbox{$\lambda^{\flat}:\T{M} \rightarrow \T^{\ast}{M}$}\, is the vector bundle morphism given by \,\mbox{$X \mapsto \ii_{X}\lambda$}.\, The bivector field $\overline{\Pi}$ is called the \emph{$\lambda$--gauge transformation} of $\Pi$ \cite{GaugeWe,GaugeBursz,Gauge}. A pair of bivector fields $\Pi$ and $\overline{\Pi}$ on $M$ are said to be \emph{gauge equivalent} if they are related by (\ref{EcGaugeEquiv}) for some differential 2--form $\lambda$ on $M$ satisfying (\ref{EcGaugeInvertible}).
If $\Pi$ is a Poisson bivector field, then $\overline{\Pi}$ is a Poisson bivector field if and only if $\lambda$ is closed along the symplectic leaves of $\Pi$. A gauge transformation modifies only the leaf-wise symplectic form of $\Pi$ by means of the pull-back of $\lambda$, preserving the characteristic foliation. Furthermore, gauge transformations preserve unimodularity.
The function \textsf{gauge\_transformation} computes the gauge transformation of a bivector field.
\begin{algorithm}[H]
    \captionsetup{justification=centering}
    \caption{\ \textsf{gauge\_transformation}(\emph{bivector, two\_form})} \label{AlgGaugeTrans}
        \rule{\textwidth}{0.4pt}
    \Input{a bivector field and a differential 2-form}
    \Output{a bivector field which is the gauge transformation induced by the differential 2-form of the given bivector field}
        \rule{\textwidth}{0.4pt}
    \begin{algorithmic}[1] 
        \Procedure{}{}
            \State $m$ $\gets$ dimension of the manifold
                \CommentNew{Given by an instance of \textsf{PoissonGeometry}}
            \State \bluecolor{bivector} $\gets$ a dictionary that represents a bivector field according to (\ref{EcMultivectorDic})
            \State \bluecolor{two\_form} $\gets$ a dictionary that represents a differential 2-form according to (\ref{EcMultivectorDic})
            \State \bluecolor{bivector\_matrix} $\gets$ \textsf{bivector\_to\_matrix}(\bluecolor{bivector})
                \CommentNew{See Algorithm \ref{AlgMatrixPoisson}}
            \State \bluecolor{2\_form\_matrix} $\gets$ \textsf{bivector\_to\_matrix}(\bluecolor{two\_form})
            \State \bluecolor{identity} $\gets$ the $m \times m$ identity matrix
            \If {det(\bluecolor{identity} $-$ \bluecolor{2\_form\_matrix} * \bluecolor{bivector\_matrix}) == 0}
                \State \textbf{return} False \CommentNew{Means that (\ref{EcGaugeInvertible}) is not invertible}
            \Else
                \State \bluecolor{gauge\_matrix} $\gets$ \bluecolor{bivector} * (\bluecolor{identity}
                    - \bluecolor{2\_form\_matrix} * \bluecolor{bivector\_matrix})
                \State \bluecolor{gauge\_bivector} $\gets$ an empty dictionary dict()
                \For {$1 \leq i < j \leq m$}
                    \State \bluecolor{gauge\_bivector}[$(i,j)$] $\gets$ \bluecolor{gauge\_matrix}[$i-1,j-1$]
                \EndFor
                \State \textbf{return} \mbox{\bluecolor{gauge\_bivector}, det(\bluecolor{identity}
                    - \bluecolor{ 2\_form\_matrix} * \bluecolor{bivector\_matrix})}
            \EndIf
        \EndProcedure
    \end{algorithmic}
\end{algorithm}
\vspace{-.5cm}
Observe that the function \textsf{gauge\_transformation} can be used to compute the gauge transformation induced by a closed differential 2-form of a Poisson bivector field.

    \subsection{Classification of Lie-Poisson bivector fields on $\R{3}$} \label{subsec:classif}
A \emph{Lie-Poisson} bivector field is a homogeneous Poisson bivector field (\ref{EcPiHomog}) for which each $\Pi^{ij}$ is a linear polynomial \cite{Konstant,Ginzburg,Dufour}. A pair of homogeneous Poisson bivector fields $\Pi$ \,and \, $\widetilde{\Pi}$ on $\R{m}$ are said to be equivalent (or isomorphic) if there exists an invertible linear operator \,$T:\R{m} \rightarrow \R{m}$\, such that
     \begin{equation}\label{EcEquivHomogRn}
          \widetilde{\Pi} \,=\, T^{\ast}\Pi.
     \end{equation}
Under this equivalence relation in the $3$-dimensional case there exist 9 non-trivial equivalence classes of Lie-Poisson bivector fields \cite{LiuXU-92}.

The function \textsf{linear\_normal\_form\_R3} computes a normal form of a Lie-Poisson bivector field on $\R{3}$. The normal forms are based on well-known classifications of (real) 3-dimensional Lie algebra isomorphisms \cite{LiuXU-92}.
\begin{algorithm}[H]
    \captionsetup{justification=centering}
    \caption{\ \textsf{linear\_normal\_form\_R3}(\emph{bivector})} \label{AlgLinNormalFormR3}
        \rule{\textwidth}{0.4pt}
    \Input{a Lie-Poisson bivector field on $\R{3}$}
    \textbf{Output:} {a linear normal form for the Lie-Poisson bivector field} \\[-0.25cm]
        \rule{\textwidth}{0.4pt}
    \begin{algorithmic}[1] 
        \Procedure{}{}
                \breakalg
            \State \bluecolor{bivector} $\gets$ a dictionary \{$(1,2)$: $\mathsf{\Pi}^{12}$, $(1,3)$: $\mathsf{\Pi}^{13}$, $(2,3)$: $\mathsf{\Pi}^{23}$\} that represents a Lie-Poisson bivector field on $\R{3}$ according to (\ref{EcMultivectorDic})
            \State \textsc{Convert} each value in \bluecolor{bivector} to symbolic expression
            \State \bluecolor{parameter} $\gets$ \textsf{x1} * $\mathsf{\Pi}^{23}$ - \textsf{x2} * $\mathsf{\Pi}^{13}$ + \textsf{x3} * $\mathsf{\Pi}^{12}$
            \State \bluecolor{hessian\_parameter} $\gets$ Hessian matrix of \bluecolor{parameter}
            \If {\textsf{modular\_vf}(\bluecolor{bivector}) == 0} \CommentNew{See Algorithm \ref{AlgModularVF}}
                \If {$\mathrm{rank}$(\bluecolor{hessian\_parameter}) == $0$}
                    \State \textbf{return} \{0: 0\} \CommentNew{A dictionary with zero key and value}
                \ElsIf {$\mathrm{rank}$(\bluecolor{hessian\_parameter}) == $1$}
                    \State \textbf{return} \{$(1,2)$: 0, $(1,3)$: 0, $(2,3)$: \textsf{x1}\}
                \ElsIf {$\mathrm{rank}$(\bluecolor{hessian\_parameter}) == $2$}
                    \If {$\mathrm{index}$(\bluecolor{hessian\_parameter}) == $2$}
                            \CommentNew{Index of quadratic forms}
                        \State \textbf{return} \{$(1,2)$: 0, $(1,3)$: \textsf{-x2}, $(2,3)$: \textsf{x1}\}
                    \Else
                        \State \textbf{return} \{$(1,2)$: 0, $(1,3)$: \textsf{x2}, $(2,3)$: \textsf{x1}\}
                    \EndIf
                \Else
                    \If {$\mathrm{index}$(\bluecolor{hessian\_parameter}) == $3$}
                            \CommentNew{Index of quadratic forms}
                        \State \textbf{return} \{$(1,2)$: \textsf{x3}, $(1,3)$: \textsf{-x2}, $(2,3)$: \textsf{x1}\}
                    \Else
                        \State \textbf{return} \{$(1,2)$: \textsf{-x3}, $(1,3)$: \textsf{-x2}, $(2,3)$: \textsf{x1}\}
                    \EndIf
                \EndIf
            \Else
                \If {$\mathrm{rank}$(\bluecolor{hessian\_parameter}) == $0$}
                    \State \textbf{return} \{$(1,2)$: 0, $(1,3)$: \textsf{x1}, $(2,3)$: \textsf{x2}\}
                \ElsIf {$\mathrm{rank}$(\bluecolor{hessian\_parameter}) == $1$}
                    \State \textbf{return} \{$(1,2)$: 0, $(1,3)$: \textsf{x1}, $(2,3)$: \textsf{4*x1 + x2}\}
                \Else
                    \If {$\mathrm{index}$(\bluecolor{hessian\_parameter}) == $2$} \CommentNew{Index of quadratic forms}
                        \State \textbf{return} \{$(1,2)$: 0, $(1,3)$: \textsf{x1 - 4*a*x2}, $(2,3)$: \textsf{4*a*x1 + x2}\}
                    \Else
                        \State \textbf{return} \{$(1,2)$: 0, $(1,3)$: \textsf{x1 + 4*a*x2}, $(2,3)$: \textsf{4*a*x1 + x2}\}
                    \EndIf
                \EndIf
            \EndIf
        \EndProcedure
    \end{algorithmic}
\end{algorithm}
\vspace{-.5cm}

    \subsection{Isomorphic Lie-Poisson Tensors on $\R{3}$} \label{subsec:isomorphic}
Using the function \textsf{isomorphic\_lie\_poisson\_R3} we can verify whether two Lie-Poisson bivector fields on $\R{3}$ are isomorphic (\ref{EcEquivHomogRn}), or not.
\begin{algorithm}[H]
    \captionsetup{justification=centering}
    \caption{\ \textsf{isomorphic\_lie\_poisson\_R3}(\emph{bivector\_1, bivector\_2})} \label{AlgIsomorphicLiePoissonR3}
        \rule{\textwidth}{0.4pt}
    \textbf{Input:} two Lie-Poisson bivector fields\\
    \textbf{Output:} verify if the Lie-Poisson bivector fields are isomorphic or not \\[-0.25cm]
        \rule{\textwidth}{0.4pt}
    \begin{algorithmic}[1] 
        \Procedure{}{}
            \State \bluecolor{bivector\_1} $\,\gets\,$ a dictionary that represents a bivector field according to (\ref{EcMultivectorDic})
                \breakalg
            \State \bluecolor{bivector\_2} $\,\gets\,$ a dictionary that represents a bivector field according to (\ref{EcMultivectorDic})
            \If {\mbox{ \textsf{linear\_normal\_form\_R3}(\bluecolor{bivector\_1})
                \,==\, \textsf{linear\_normal\_form\_R3}(\bluecolor{ bivector\_2})}}
                    \CommentNew{See Algorithm \ref{AlgLinNormalFormR3}}
                \State \textbf{return} True
            \Else
                \State \textbf{return} False
            \EndIf
        \EndProcedure
    \end{algorithmic}
\end{algorithm}
\vspace{-.5cm}

    \subsection{Flaschka-Ratiu Bivector Fields} \label{sec:flaschka}
Given $m-2$ functions \,$K_{1},...,K_{m-2} \in \Cinf{M}$ on an oriented $m$-dimensional manifold $(M,\Omega)$, with volume form $\Omega$, we can construct a Poisson bivector field $\Pi$ on $M$ defined by
    \begin{equation*}
        \ii_{\Pi}\Omega \,:=\, \dd{K_{1}} \wedge \cdots \wedge \dd{K_{m-2}}.
    \end{equation*}
Clearly, $\Pi$ is non-trivial on the open subset of $M$ where \,$K_{1},\ldots,K_{m-2}$\, are (functionally) independent. Moreover, by construction, each $K_{l}$ is a Casimir function of $\Pi$. These class of Poisson bivector fields are called \emph{Flaschka-Ratiu} bivector fields \cite{Damianou}. In coordinates, if \,$\Omega=\dd{x^{1}} \wedge \cdots \wedge \dd{x^{m}}$,\, then
    \begin{equation}\label{EcFlaschkaFormula}
        \Pi \,=\, (-1)^{i+j}\,\det{P_{[i,j]}}\cdot\frac{\partial}{\partial{x^{i}}} \wedge \frac{\partial}{\partial{x^{j}}}, \qquad
            1 \leq i < j \leq m.
    \end{equation}
 Here $P$ denotes the \,$(m-2) \times m$-matrix whose $k$-th row is \,$\big(\partial{K_{k}}/\partial{x^1},\ldots,\partial{K_{k}}/\partial{x^m}\big)$,\, for \,$k=1,\ldots,m-2$;\, and $P_{[i,j]}$ the matrix $P$ without the columns $i$ and $j$. Moreover, the symplectic form $\omega_{S}$ of $\Pi$ on a 2-dimensional (symplectic) leaf \,$S \subseteq M$\, is given by
    \begin{equation}\label{EcSymplecticForm}
        \omega_{S} \,=\, \tfrac{1\ }{|\Pi|^{2}}\left[\, (-1)^{i+j+1}\det{P_{[i,j]}}\cdot\dd{x^{i}} \wedge \dd{x^{j}}\,\right]\Big|_{S}, \quad |\Pi|^{2} := \sum_{1 \leq i < j \leq m}\big(\det{P_{[i,j]}}\,\big)^{2}.
    \end{equation}

The function \textsf{flaschka\_ratiu\_bivector} computes the Flaschka-Ratiu bivector field and the corresponding symplectic form of a `maximal' set of scalar functions \cite{Damianou, Nar-2015, PabloWrinFib, PabSua-2018}.
\begin{algorithm}[H]
    \captionsetup{justification=centering}
    \caption{\ \textsf{flaschka\_ratiu\_bivector}(\emph{casimir\_list})} \label{AlgFlaschka}
        \rule{\textwidth}{0.4pt}
    \Input{$m-2$\, scalar functions}
    \Output{the Flaschka-Ratiu bivector field induced by the $m-2$ functions and the symplectic form of this Poisson bivector field}
        \rule{\textwidth}{0.4pt}
    \begin{algorithmic}[1] 
        \Procedure{}{}
            \State $m$ $\gets$ dimension of the manifold
                \CommentNew{Given by an instance of \textsf{PoissonGeometry}}
            \State \bluecolor{casimir\_list} $\gets$ a list \,\textsf{[`K1', \ldots, `K\{$m-2$\}']}\, with $m-2$ string expressions
            \Statex \CommentNew{Each string expression represents a scalar function}
            \If {at least two functions in \bluecolor{casimir\_list} are functionally dependent}
                \State \textbf{return} \{0: 0\} \CommentNew{A dictionary with zero key and value}
            \Else
                \State \bluecolor{matrix\_gradients} $\gets$ a symbolic $(m-2) \times m$--matrix
                \For {$i=1$ to $m-2$}
                    \State \textsc{Convert} \textsf{K}{$i$} to symbolic expressions
                    \State \textsc{Compute} the gradient vector $\nabla$\textsf{K}{$i$} of \textsf{K}{$i$}
                        \breakalg
                    \State \textsc{Append} $\nabla$\textsf{K}{$i$} to \bluecolor{matrix\_gradients} as its $i$-th row
                \EndFor
                \State \bluecolor{flaschka\_bivector} $\gets$ an empty dictionaty dict()
                \State \bluecolor{sum\_bivector} $\gets$ 0
                \For {each \,$1 \leq i < j \leq m$\,}
                    \State \textsc{Remove} from \bluecolor{matrix\_gradients} the $i$-th and $j$-th columns
                    \State \bluecolor{flaschka\_bivector}[$(i,j)$] $\gets$ $(-1)^{i+j}$ * det(\bluecolor{matrix\_gradients})
                        \CommentNew{See (\ref{EcFlaschkaFormula})}
                    \State \bluecolor{sum\_bivector} $\gets$ \bluecolor{sum\_bivector} + det(\bluecolor{matrix\_gradients})**2
                    \State \textsc{Append} to \bluecolor{matrix\_gradients} the $i$-th and $j$-th removed
                \EndFor
                \State \bluecolor{symplectic\_form} $\gets$ an empty dictionary dict()
                \For {each \,$1 \leq i < j \leq m$\,}
                    \State \bluecolor{symplectic\_form}[$(i,j)$] $\gets$ $(-1)$ * \bluecolor{sum\_bivector}
                        * \bluecolor{ flaschka\_bivector}[$(i,j)$]
                    \Statex \CommentNew{See (\ref{EcSymplecticForm})}
                \EndFor
                \State \textbf{return} \bluecolor{flaschka\_bivector}, \bluecolor{symplectic\_form}
            \EndIf
        \EndProcedure
    \end{algorithmic}
\end{algorithm}
\vspace{-.5cm}

    \subsection{Test Type Functions} \label{sec:test}
In this section we describe our implementation of some useful functions in the \textsf{PoissonGeometry} module which allow us to verify whether a given geometric object on a Poisson manifold satisfies certain property. The algorithms for each of these functions are similar, as they are decision-making processes.

    \subsubsection{Jacobi Identity} \label{subsec:jacobidentity}
We can verify in \textsf{PoissonGeometry} if a given bivector field $\Pi$ is a Poisson bivector field or not.
\begin{algorithm}[H]
    \captionsetup{justification=centering}
    \caption{\ \textsf{is\_poisson\_tensor}(\emph{bivector})} \label{AlgIsPoissonTensor}
        \rule{\textwidth}{0.4pt}
    \textbf{Input:} a bivector field \\
    \textbf{Output:} verify if the bivector field is a Poisson bivector field or not \\[-0.25cm]
        \rule{\textwidth}{0.4pt}
    \begin{algorithmic}[1] 
        \Procedure{}{}
            \State \bluecolor{bivector} $\,\gets\,$ a dictionary that represents a bivector field according to (\ref{EcMultivectorDic})
            \If {\textsf{lichnerowicz\_poisson\_operator}(\bluecolor{bivector}, \bluecolor{bivector}) == \textsf{\{0: 0\}}}
                \Statex \CommentNew{See Algorithm \ref{AlgDeltaPi}}
                \State \textbf{return} True
            \Else
                \State \textbf{return} False
            \EndIf
        \EndProcedure
    \end{algorithmic}
\end{algorithm}
\vspace{-.5cm}

    \subsubsection{Kernel of a Bivector Field} \label{subsec:kernel}
The kernel of a bivector field $\Pi$ is the subspace \,\mbox{$\ker{\Pi}:=\{\alpha \in \T^{\ast}M \,|\, \Pi^{\natural}(\alpha) = 0\}$}\, of $\T^{\ast}M$. It is defined as the kernel of its sharp morphism (\ref{EcPiSharp}), and is defined likewise  for Poisson bivector fields \cite{Dufour, Camille}.


\begin{algorithm}[H]
    \captionsetup{justification=centering}
    \caption{\ \textsf{is\_in\_kernel}(\emph{bivector, one\_form})} \label{AlgIsKernel}
        \rule{\textwidth}{0.4pt}
    \textbf{Input:} a bivector field and a diferential 1--form \\
    \Output{verify if the differential 1--form belongs to the kernel of the (Poisson) bivector field}
        \rule{\textwidth}{0.4pt}
    \begin{algorithmic}[1] 
        \Procedure{}{} \vspace{2pt}
            \State \bluecolor{bivector} $\,\gets\,$ a dictionary that represents a bivector field according to (\ref{EcMultivectorDic})
            \State \bluecolor{one\_form} $\,\gets\,$ a dictionary that represents a differential 1-form according to (\ref{EcMultivectorDic})
            \If {\textsf{sharp\_morphism}(\bluecolor{bivector}, \bluecolor{one\_form}) == \textsf{\{0: 0\}}}
                    \CommentNew{See Algorithm \ref{AlgSharpMorp}}
                \State \textbf{return} True
            \Else
                \State \textbf{return} False
            \EndIf
        \EndProcedure
    \end{algorithmic}
\end{algorithm}
\vspace{-.5cm}

    \subsubsection{Casimir Functions} \label{subsec:casimir}
A function \,\mbox{$K \in \Cinf{M}$}\, is said to be a Casimir funtion of a Poisson bivector field $\Pi$ if its Hamiltonian vector field (\ref{EcPiHamVectField}) is zero. Equivalently, if its exterior derivative $\dd{K}$ belongs to the kernel of $\Pi$ \cite{Dufour,Damianou,Camille}. 
\begin{algorithm}[H]
    \captionsetup{justification=centering}
    \caption{\ \textsf{is\_casimir}(\emph{bivector, function})} \label{AlgIsCasimir}
        \rule{\textwidth}{0.4pt}
    \textbf{Input:} a Poisson bivector field and a scalar function \\
    \Output{verify if the scalar function is a Casimir function of the Poisson bivector field}
        \rule{\textwidth}{0.4pt}
    \begin{algorithmic}[1] 
        \Procedure{}{}
            \State \bluecolor{bivector} $\,\gets\,$ a dictionary that represents a bivector field according to (\ref{EcMultivectorDic})
            \State \bluecolor{function} $\,\gets\,$ a string expression
                \CommentNew{Represent a scalar function}
            \If {\textsf{hamiltonian\_vf}(\bluecolor{bivector}, \bluecolor{function}) == \textsf{\{0: 0\}}}
                \Statex \CommentNew{See Algorithm \ref{AlgHamVectField}}
                \State \textbf{return} True
            \Else
                \State \textbf{return} False
            \EndIf
        \EndProcedure
    \end{algorithmic}
\end{algorithm}
\vspace{-.5cm}

    \subsubsection{Poisson Vector Fields} \label{subsec:poissonvf}
A vector field $W$ on $M$ is said to be a Poisson vector field of a Poisson bivector field $\Pi$ if it commutes with respect to the Schouten-Nijenhuis bracket, \,$\cSch{W,\Pi} \,=\, 0$ \cite{Dufour,Camille}.
\begin{algorithm}[H]
    \captionsetup{justification=centering}
    \caption{\ \textsf{is\_poisson\_vf}(\emph{bivector, vector\_field})} \label{AlgIsPoissonVectorField}
        \rule{\textwidth}{0.4pt}
    \textbf{Input:} a Poisson bivector field and a vector field \\
    \textbf{Output:} {verify if the vector field is a Poisson vector field of the Poisson bivector field} \\[-0.25cm]
        \rule{\textwidth}{0.4pt}
\begin{algorithmic}[1] 
    \Procedure{}{}
            \State \bluecolor{bivector} $\,\gets\,$ a dictionary that represents a bivector field according to (\ref{EcMultivectorDic})
            \State \bluecolor{vector\_field} $\,\gets\,$ a dictionary that represents a vector field according to (\ref{EcMultivectorDic})
                \breakalg
            \If {\textsf{lichnerowicz\_poisson\_operator}(\bluecolor{bivector}, \bluecolor{vector\_field}) == \textsf{\{0: 0\}}}
                \Statex \CommentNew{See Algorithm \ref{AlgDeltaPi}}
                \State \textbf{return} True
            \Else
                \State \textbf{return} False
            \EndIf
\EndProcedure
    \end{algorithmic}
\end{algorithm}
\vspace{-.5cm}

    \subsubsection{Poisson Pairs} \label{subsec:poissonpair}
We can verify whether a couple of Poisson bivector fields $\Pi$ and $\Psi$ form a Poisson pair. That is, if the sum $\Pi + \Psi$ is again a Poisson bivector field or, equivalently, if $\Pi$ and $\Psi$ commute with respect to the Schouten-Nijenhuis bracket, \,$\cSch{\Pi, \Psi} \,=\, 0$\, \cite{Dufour,Camille}.
\begin{algorithm}[H]
    \captionsetup{justification=centering}
    \caption{\ \textsf{is\_poisson\_pair}(\emph{bivector\_1, bivector\_2})} \label{AlgIsPoissonPair}
        \rule{\textwidth}{0.4pt}
    \textbf{Input:} two Poisson bivector fields. \\
    \Output{verify if the bivector fields commute with respect to the Schouten-Nijenhuis bracket}
        \rule{\textwidth}{0.4pt}
\begin{algorithmic}[1] 
    \Procedure{}{}
            \State \bluecolor{bivector\_1} $\,\gets\,$ a dictionary that represents a bivector field according to (\ref{EcMultivectorDic})
            \State \bluecolor{bivector\_2} $\,\gets\,$ a dictionary that represents a bivector field according to (\ref{EcMultivectorDic})
            \If {\textsf{lichnerowicz\_poisson\_operator}(\bluecolor{bivector\_1}, \bluecolor{bivector\_2}) == \textsf{\{0: 0\}}}
                \Statex \CommentNew{See Algorithm \ref{AlgDeltaPi}}
                \State \textbf{return} True
            \Else
                \State \textbf{return} False
            \EndIf
\EndProcedure
    \end{algorithmic}
\end{algorithm}
\vspace{-.5cm}

\section{\textsf{PoissonGeometry}: Syntax and Applications} \label{sec:class}
\textsf{PoissonGeometry} is our python module for local calculus on Poisson manifolds. First we define a tuple of symbolic variables that emulate local coordinates on a finite (Poisson) smooth manifold $M$. By default, these symbolic variables are just the juxtaposition of the symbol \textsf{x} and an index of the set $\{1,\ldots,\mbox{$m=\dim{M}$}\}$: (\text{\textsf{x1}}, \ldots, \text{\textsf{xm}}).

\noindent\textbf{Scalar Functions.} A local representation of a scalar function in \textsf{PoissonGeometry} is written using \emph{string literal expressions}. For example, the function \,$f = a(x^1)^2 + b(x^2)^2 + c(x^3)^2$\, should be written exactly as follows: \textsf{`a * x1**2 + b * x2**2 + c * x3**2'}. It is \textsl{important} to remember that all characters that are not local coordinates are treated as (symbolic) parameters: \textsf{a}, \textsf{b} and \textsf{c} for the previous example.

\noindent\textbf{Multivector Fields and Differential forms.} Both multivector fields and differential forms are written using \emph{dictionaries} with \textsl{tuples of integers} as \textsl{keys} and \textsl{string} type \textsl{values}. If the coordinate expression of an $a$--multivector field $A$ on $M$, with \,$a \in \mathbb{N}$,\, is given by,
    \begin{equation*}
        A \,= \sum_{1 \leq i_1 < i_2 < \cdots < i_a \leq m} A^{i_1 i_2 \cdots i_a}\,\frac{\partial}{\partial{x^{i_1}}} \wedge \frac{\partial}{\partial{x^{i_2}}} \wedge \cdots \wedge \frac{\partial}{\partial{x^{i_a}}}, \quad A^{i_1 \cdots i_a} \,=\, A^{i_1 \cdots i_a}(x),
    \end{equation*}
then $A$ should be written using a dictionary, as follows:
    \begin{equation}\label{EcMultivectorDic}
        \Big\{(1,...,a): \mathscr{A}^{1\cdots a}, \,...,\, (i_1,...,i_a): \mathscr{A}^{i_1 \cdots i_a}, \,...,\, (m-a+1,...,m): \mathscr{A}^{m-a+1\cdots m}\Big\}.
    \end{equation}   
Here each key $(i_1,\ldots,i_a)$ is a tuple containing ordered indices \mbox{$1 \leq i_1 < \cdots < i_a \leq m$} and the corresponding value $\mathscr{A}^{i_1 \cdots i_a}$ is the string expression of the scalar function (coefficient) $A^{i_1 \cdots i_a}$ of $A$.

The syntax for differential forms is the same. It is important to remark that we can only write the keys and values of \textsl{non-zero coefficients}. See the documentation for more details.

    \subsection{Applications}
We will now describe two applications. One of \textsf{gauge\_trans-formation}, used here to derive a characterization of gauge transformations on $\R{3}$ (see, Subsection \ref{subsec:gauge}), and a second one of \textsf{jacobiator}, used here to construct a family of Poisson bivector fields on $\R{4}$ (\ref{EcJacobiPi}).

\noindent\textbf{Gauge Transformations on $\R{3}$.} For an arbitrary bivector field on $\R{3}_{x}$,
    \begin{equation}\label{EcPiR3}
        \Pi \,=\, \Pi^{12}\,\frac{\partial}{\partial{x^{1}}} \wedge \frac{\partial}{\partial{x^{2}}}
            \,+\, \Pi^{13}\,\frac{\partial}{\partial{x^{1}}} \wedge \frac{\partial}{\partial{x^{3}}}
            \,+\, \Pi^{23}\,\frac{\partial}{\partial{x^{2}}} \wedge \frac{\partial}{\partial{x^{3}}},
    \end{equation}
and an arbitrary differential 2-form
    \begin{equation}\label{EcLambdaR3}
        \lambda \,=\, \lambda_{12}\,\dd{x^{1}} \wedge \dd{x^{2}} \,+\, \lambda_{13}\,\dd{x^{1}} \wedge \dd{x^{3}}
                \,+\, \lambda_{23}\,\dd{x^{2}} \wedge \dd{x^{3}},
    \end{equation}
we compute:
\begin{tcolorbox}[arc=0mm, boxsep=0mm, skin=bicolor, colback=blue!20, colframe=blue!25, colbacklower=blue!1]
$>>>$ \textsf{pg3 = PoissonGeometry(3)} \\
$>>>$ \textsf{P = \{(1,2): `P12', (1,3): `P13', (2,3): `P23'\}} \\
$>>>$ \textsf{lambda = \{(1,2): `L12', (1,3): `L13', (2,3): `L23'\}} \\
$>>>$ \textsf{(gauge\_bivector, determinant) = pg3.gauge\_transformation(P, lambda)} \\
$>>>$ \textsf{print(gauge\_bivector)} \\
$>>>$ \textsf{print(determinant)}
    \tcblower
\textsf{\big\{}\parbox[t]{0.85\linewidth}{\textsf{(1,2): P12/(L12*P12 + L13*P13 + L23*P23 + 1), \\
(1,3): P13/(L12*P12 + L13*P13 + L23*P23 + 1), \\
(2,3): P23/(L12*P12 + L13*P13 + L23*P23 + 1)\big\}}} \\[0.15cm]
\textsf{(L12*P12 + L13*P13 + L23*P23 + 1)**2}
\end{tcolorbox}
The symbols \textsf{P12}, \textsf{P13}, \textsf{P23}, and \textsf{L12}, \textsf{L13}, \textsf{L23} stand for the coefficients of $\Pi$ and $\lambda$, in that order. Then, (see (\ref{EcGaugeInvertible})):
    \begin{equation}\label{EcDeterminantR3}
        \det\big(\mathrm{Id} - \lambda^{\flat} \circ \Pi^{\natural}\big) \,=\, \big(\lambda_{12}\Pi^{12} + \lambda_{13}\Pi^{13} + \lambda_{23}\Pi^{23} + 1\big)^2
    \end{equation}
So, for \,$1 \leq i < j \leq 3$, the $\lambda$-gauge transformation $\overline{\Pi}$ of $\Pi$ is given by:
    \begin{equation}\label{EcPiGaugeR3}
        \overline{\Pi} \,=\, \frac{\Pi^{12}}{\lambda_{ij}\Pi^{ij} + 1}\,\frac{\partial}{\partial{x^{1}}} \wedge \frac{\partial}{\partial{x^{2}}}
            \,+\, \frac{\Pi^{13}}{\lambda_{ij}\Pi^{ij} + 1}\,\frac{\partial}{\partial{x^{1}}} \wedge \frac{\partial}{\partial{x^{3}}} \,+\, \frac{\Pi^{23}}{\lambda_{ij}\Pi^{ij} + 1}\,\frac{\partial}{\partial{x^{2}}} \wedge \frac{\partial}{\partial{x^{3}}}
    \end{equation}

With these ingredients, we can now show:

\begin{proposition}\label{Prop:gauge}
Let $\Pi$ be a bivector field on a 3--dimensional smooth manifold $M$. Then, given a differential 2--form $\lambda$ on $M$, the $\lambda$-gauge transformation $\overline{\Pi}$ (\ref{EcGaugeEquiv}) of $\Pi$ is well defined on the open subset,
    \begin{equation}\label{EcFdetGauge}
        \big\{ F := \big\langle \lambda,\Pi \big\rangle + 1 \,\neq\, 0 \big\} \,\subseteq\, M.
    \end{equation} 
Moreover, $\overline{\Pi}$ is given by
    \begin{equation}\label{EcPiGauge3}
        \overline{\Pi} \,=\, \tfrac{1}{F}\Pi.
    \end{equation}
If $\Pi$ is Poisson and $\lambda$ is closed along the leaves of $\Pi$, then $\overline{\Pi}$ is also Poisson.
\end{proposition}

\begin{proof}
Suppose (\ref{EcPiR3}) and (\ref{EcLambdaR3}) are coordinate expressions of $\Pi$ and $\lambda$ on a chart $(U;x^{1},x^{2},x^{3})$ of $M$. Observe that the pairing of $\Pi$ and $\lambda$ is given by 
\[
\langle \lambda,\Pi \rangle = \lambda_{12}\Pi^{12} + \lambda_{13}\Pi^{13} + \lambda_{23}\Pi^{23}.
\]
Hence (\ref{EcDeterminantR3}) yields \,$\det(\mathrm{Id} - \lambda^{\flat} \circ \Pi^{\natural}) = (\langle \lambda,\Pi \rangle + 1)^2$.\, This implies that the morphism $\mathrm{Id} - \lambda^{\flat} \circ \Pi^{\natural}$ is invertible on the open subset in (\ref{EcFdetGauge}) and, in consequence, the $\lambda$--gauge transformation of $\Pi$ (\ref{EcGaugeInvertible}). Finally, formula (\ref{EcPiGauge3}) follows from (\ref{EcPiGaugeR3}) as
\[
\overline{\Pi} = \frac{\Pi^{12}}{F}\frac{\partial}{\partial{x^{1}}} \wedge \frac{\partial}{\partial{x^{2}}} + \frac{\Pi^{13}}{F}\frac{\partial}{\partial{x^{1}}} \wedge \frac{\partial}{\partial{x^{3}}} + \frac{\Pi^{23}}{F}\frac{\partial}{\partial{x^{2}}} \wedge \frac{\partial}{\partial{x^{3}}},
\]
for $F$ in (\ref{EcFdetGauge}.
\end{proof}

\noindent\textbf{Parametrized Poisson Bivector Fields.} Poisson bivector fields also play an important role in the theory of deformation quantization, which is linked to quantum mechanics \cite{Bayen}.
They appear in star products, that is, in general deformations of the associative algebra of smooth functions of a symplectic manifold \cite{Flato}. Our module \textsf{PoissonGeometry} can be used to study particular problems around deformations of Poisson bivector fields and star products.

For example, we can modify the following $4$--parametric bivector field on $\R{4}$
    \begin{equation*}
        \Pi \,=\, a_{1}x^2\,\frac{\partial}{\partial{x^{1}}} \wedge \frac{\partial}{\partial{x^{2}}}
            + a_{2}x^3\,\frac{\partial}{\partial{x^{1}}} \wedge \frac{\partial}{\partial{x^{3}}}
            + a_{3}x^4\,\frac{\partial}{\partial{x^{1}}} \wedge \frac{\partial}{\partial{x^{4}}}
            + a_{4}x^1\,\frac{\partial}{\partial{x^{2}}} \wedge \frac{\partial}{\partial{x^{3}}},
    \end{equation*}
using the \textsf{jacobiator} function to construct a family of Poisson bivector fields on $\R{4}$:
\begin{tcolorbox}[arc=0mm, boxsep=0mm, skin=bicolor, colback=blue!20, colframe=blue!25, colbacklower=blue!1, breakable]
$>>>$ \textsf{pg4 = PoissonGeometry(4)} \\
$>>>$ \textsf{P = \big\{(1,2): `a1*x2', (1,3): `a2*x3', (1,4): `a3*x4', (2,3): `a4*x1'\big\}} \\
$>>>$ \textsf{pg4.jacobiator(P)}
    \tcblower
\textsf{\big\{(1,2,3): -2*a4*x1*(a1 + a2), (2,3,4): -2*a3*a4*x4\big\}}
\end{tcolorbox}
Therefore
    \begin{equation*}
        \cSch{\Pi,\Pi} \,=\, -2a_{4}(a_{1} + a_{2})\,x^1\,\frac{\partial}{\partial{x^{1}}} \wedge \frac{\partial}{\partial{x^{2}}} \wedge \frac{\partial}{\partial{x^{3}}}
            \,-\, 2a_{3}a_{4}\,x^4\,\frac{\partial}{\partial{x^{2}}} \wedge \frac{\partial}{\partial{x^{3}}} \wedge \frac{\partial}{\partial{x^{4}}}
    \end{equation*}
Hence, we have two cases, explained in the following:

\begin{lemma}\label{example}If \,\mbox{$a_{4}=0$},\, then $\Pi$ determines a 3--parametric family of Poisson bivector fields on $\R{4}_{x}$:
    \begin{equation}
        \Pi \,=\, a_{1}x^2\,\frac{\partial}{\partial{x^{1}}} \wedge \frac{\partial}{\partial{x^{2}}}
            + a_{2}x^3\,\frac{\partial}{\partial{x^{1}}} \wedge \frac{\partial}{\partial{x^{3}}}
            + a_{3}x^4\,\frac{\partial}{\partial{x^{1}}} \wedge \frac{\partial}{\partial{x^{4}}}.
    \end{equation}
If \,\mbox{$a_{2}=-a_{1}$}\, and \,\mbox{$a_{3}=0$},\, then $\Pi$ determines a 2--parametric family of Poisson bivector fields on $\R{4}_{x}$:
\begin{equation}
        \Pi \,=\, a_{1}x^2\,\frac{\partial}{\partial{x^{1}}} \wedge \frac{\partial}{\partial{x^{2}}} - a_{1}x^3\,\frac{\partial}{\partial{x^{1}}} \wedge \frac{\partial}{\partial{x^{3}}} + a_{4}x^1\,\frac{\partial}{\partial{x^{2}}} \wedge \frac{\partial}{\partial{x^{3}}}.
    \end{equation}
\end{lemma}

\noindent\textbf{Related Work.} Notable contributions in similar directions include; computations of normal forms in Hamiltonian  dynamics (in Maxima) \cite{TV-19}, symbolic tests of the Jacobi identity for generalized Poisson brackets  and their relation to hydrodynamics \cite{JacobiMath}, and an implementation of the Schouten-Bracket for multivector fields (in Sage \footnote{\tt  https://trac.sagemath.org/ticket/23429 }). 

Our work here is, to the best of our knowledge, the first comprehensive implementation of routine computations used in Poisson geometry, and in Python (based on \textsf{SymPy}\cite{Sympy}).

\noindent\textbf{Future Directions.} With the algorithms in this paper, numerical extensions for the same methods can be developed. Explicit computations of Poisson cohomology can also be explored. These are the subjects of ongoing, and forthcoming work.

    \section*{Acknowledgments}
We thank Luis A. G\'omez-Telesforo for many fruitful discussions.

\bibliographystyle{siamplain}

\end{document}

%% file: ex_shared.tex
\usepackage{amsmath}
\usepackage{lipsum}
\usepackage{amsfonts}
\usepackage{graphicx}
\usepackage{float}
\usepackage{caption}
\usepackage{epstopdf}
\usepackage{algpseudocode}
\usepackage{algorithmicx}
\ifpdf
  \DeclareGraphicsExtensions{.eps,.pdf,.png,.jpg}
\else
  \DeclareGraphicsExtensions{.eps}
\fi

\newsiamremark{remark}{Remark}
\newsiamremark{hypothesis}{Hypothesis}
\newsiamremark{example}{Example}
\crefname{hypothesis}{Hypothesis}{Hypotheses}
\newsiamthm{claim}{Claim}

\usepackage[mathscr]{eucal}
\usepackage{sansmath}
\usepackage{verbatim}
\usepackage{tcolorbox}
\tcbuselibrary{skins}
\tcbuselibrary{breakable}
\usepackage{color}
\definecolor{bluee}{rgb}{0,0,0.5}
\definecolor{redd}{rgb}{0.6,0,0}
\definecolor{greenn}{rgb}{0,0.5,0}
\usepackage{tikz}
\usetikzlibrary{arrows}

\edef\svtheparindent{\the\parindent}
\usepackage{parskip}
\parindent=\svtheparindent\relax

\newcommand{\R}[1]{\mathbb{R}^{#1}}
\newcommand{\T}{\mathsf{T}}
\newcommand{\dd}{\mathrm{d}}
\newcommand{\ii}{\mathbf{i}}

\newcommand{\cSch}[1]{[\hspace{-0.065cm}[ #1 ]\hspace{-0.065cm}]}
\newcommand{\Cinf}[1]{\mathbf{\mathit{C}}^{\infty}_{#1}}
\newcommand{\bluecolor}[1]{{\color{blue} #1}}
\newcommand{\Input}[1]{\textbf{Input:\ }\parbox[t]{0.91\linewidth}{#1}}
\newcommand{\Output}[1]{\textbf{Output:\ }\parbox[t]{0.88\linewidth}{#1\vspace{2pt}}}
\newcommand{\CommentNew}[1]{\Comment{\mbox{\scriptsize\color{greenn}{#1}}}} 
\newcommand{\breakalg}{\algstore{bkbreak}\end{algorithmic}\end{algorithm}\begin{algorithm}[H]\vspace{0.1cm}\begin{algorithmic}[1]\algrestore{bkbreak}}
\DeclareSymbolFont{sfletters}{OML}{cmbrm}{m}{it}
\DeclareMathSymbol{\salpha}{\mathord}{sfletters}{"0B}
\DeclareMathSymbol{\sbeta}{\mathord}{sfletters}{"0C}

\headers{Computational Poisson Geometry I}{Evangelista-Alvarado, Ru\'iz-Pantale\'on, and Su\'arez-Serrato}

\title{On Computational Poisson geometry I: \\ Symbolic Foundations
            \thanks{
            This research was partially supported by CONACyT, ``Programa para un Avance Global e Integrado de la Matemática Mexicana'' FORDECYT 265667 and UNAM-DGAPA-PAPIIT-IN104819. }}

\author{M. A. Evangelista-Alvarado$^{1 \dag}$,
        J. C. Ru\'iz-Pantale\'on$^{2 \dag}$, \\  
        and P. Su\'arez-Serrato$^{3}$
    \thanks{Instituto de Matem\'aticas, Universidad Nacional Aut\'onoma de M\'exico (UNAM), M\'exico.
            $^{1}$~\email{\mbox{miguel.eva.alv\,@matem.unam.mx}},
            $^{2}$~\email{\mbox{jcpanta\,@im.unam.mx}},
            $^{3}$~\email{\mbox{pablo\,@im.unam.mx}}.}}

\usepackage{amsopn}
